\documentclass[oneside,english,leqno]{amsart}
\usepackage[T1]{fontenc}
\usepackage[latin1]{inputenc}
\usepackage{amssymb}

\makeatletter
 \theoremstyle{plain}
\newtheorem{thm}{Theorem}[section]

\usepackage{babel}
\makeatother
\begin{document}

\title{An Inequality for Basic Confluent Hypergeometric Series}

\author{Ruiming Zhang}

\curraddr{School of Mathematics\\
Guangxi Normal University\\
Guilin City, Guangxi 541004\\
P. R. China.}

\email{ruimingzhang@yahoo.com}

\begin{abstract}
In this very short note we will derive an inequality for a class of
entire functions including all the confluent basic hypergeometric
series and an inequality for a class of meromorphic functions including
theta functions. 
\end{abstract}

\subjclass{\noindent Primary 30E15. Secondary 33D45. }

\keywords{\noindent Basic confluent hypergeometric series, $q$-Airy function,
Ramanujan's entire function, inequalities, entire functions, meromorphic
functions, theta functions. }

\maketitle

\section{Introduction\label{sec:Introduction}}

In the studies of Plancherel-Rotach type asymptotics of $q$-orthogonal
polynomials such as Ismail-Masson $\left\{ h_{n}(x|q)\right\} _{n=0}^{\infty}$,
Stieltjes-Wigert $\left\{ S_{n}(x;q)\right\} _{n=0}^{\infty}$ and
$q$-Laguerre $\left\{ L_{n}^{(\alpha)}(x;q)\right\} _{n=0}^{\infty}$\cite{Ismail1,Ismail6,Ismail7},
the Ramanujan function $A_{q}(z)$, which is also called $q$-Airy
function in the literature, plays a very important role. It is called
$q$-Airy function because it appears repeatedly in the Plancherel-Rotach
type asymptotics of $q$-orthogonal polynomials, just like classical
Airy function in the classical Plancherel-Rotach asymptotics of classical
orthogonal polynomials\cite{Ismail2,Szego}. In \cite{Zhang}, we
have investigated certain properties of the Ramanujan function $A_{q}(z)$.
It turns out that, the polynomials like Ismail-Masson $\left\{ h_{n}(x|q)\right\} _{n=0}^{\infty}$,
Stieltjes-Wigert $\left\{ S_{n}(x;q)\right\} _{n=0}^{\infty}$, $q$-Laguerre
$\left\{ L_{n}^{(\alpha)}(x;q)\right\} _{n=0}^{\infty}$ and the Ramanujan
function $A_{q}(z)$, they are all in a particular class of entire
functions including all the confluent basic hypergeometric series.
In this short note will show that this class of entire functions satisfy
an inequality similar to the inequality for the Ramanujan function
$A_{q}(z)$ in \cite{Zhang}. Finally, we will state an inequality
for a class of meromorphic functions including $\Theta(z|q)$ defined
as in \cite{Ismail6,Ismail7}. $\Theta(z|q)$ also appears in the
main terms of the Plancherel-Rotach type asymptotics when a scaling
parameter is irrational.

As in papers\cite{Ismail6,Ismail7,Zhang}, we will follow the usual
notations from $q$-series \cite{Andrews,Gasper,Ismail2}

\begin{equation}
(a;q)_{0}:=1\quad(a;q)_{n}:=\prod_{k=0}^{n}(1-aq^{k}),\quad\left[\begin{array}{c}
n\\
k\end{array}\right]_{q}:=\frac{(q;q)_{n}}{(q;q)_{k}(q;q)_{n-k}}.\label{eq:1.1}\end{equation}
Moreover, we shall use the multishifted factorials \cite{Andrews,Gasper,Ismail2}
\begin{equation}
(a_{1},a_{2},\dotsc,a_{k};q)_{n}:=\prod_{j=1}^{k}(a_{j};q)_{n}.\label{eq:1.2}\end{equation}
Thoughout this paper, we shall always assume that $0<q<1$, hence
$n=\infty$ is allowed in the above definitions.

A confluent basic hypergeometric series is formally defined as \cite{Andrews,Gasper,Ismail2}\begin{equation}
_{r}\phi_{s}\left(\begin{array}{c|c}
\begin{array}{c}
a_{1},\dotsc,a_{r}\\
b_{1},\dots,b_{s}\end{array} & q,z\end{array}\right)=\sum_{k=0}^{\infty}\frac{(a_{1},\dotsc,a_{r};q)_{k}}{(b_{1},\dotsc,b_{s};q)_{k}}z^{k}\left(-q^{(k-1)/2}\right)^{k(s+1-r)}\label{eq:1.3}\end{equation}
 with $s+1-r>0$ .

We will use a simple summation formula\begin{equation}
\frac{1}{(z;q)_{\infty}}=\sum_{k=0}^{\infty}\frac{z^{k}}{(q;q)_{k}},\label{eq:1.4}\end{equation}
which is a special case of the $q$-binomial theorem \cite{Andrews,Gasper,Ismail2},
\begin{equation}
\frac{(az;q)_{\infty}}{(z;q)_{\infty}}=\sum_{k=0}^{\infty}\frac{(a;q)_{k}}{(q;q)_{k}}z^{k}.\label{eq:1.5}\end{equation}
The Ramanujan function $A_{q}(z)$ is defined as\cite{Ismail2} \begin{equation}
A_{q}(z):=\sum_{k=0}^{\infty}\frac{q^{k^{2}}}{(q;q)_{k}}(-z)^{k}.\label{eq:1.6}\end{equation}
In \cite{Zhang} we proved that\begin{equation}
\left|A_{q}(z)\right|\le\frac{\left(\frac{|z|}{\sqrt{q}}\right)^{1/2}\exp\left\{ -\frac{\log^{2}\left|z\right|}{4\log q}\right\} }{(q;q)_{\infty}}\label{eq:1.7}\end{equation}
 for any nonzero complex number $z$ and \begin{equation}
\left|A_{q}(z)\right|\le e^{q|z|/(1-q)}\label{eq:1.8}\end{equation}
 for any complex number $z$. 

The theta function $\Theta(z|q)$ is defined as\begin{equation}
\Theta(z|q)=\sum_{k=-\infty}^{\infty}q^{k^{2}}z^{k},\label{eq:1.9}\end{equation}
 for $0<q<1$ and all nonzero complex number $z$.

\section{The Inequality\label{The Inequality}}

Let \begin{equation}
f(z)=\sum_{k=0}^{\infty}\frac{(a_{1},\dotsc,a_{r};q)_{k}q^{lk^{2}}}{(b_{1},\dotsc,b_{s},q;q)_{k}}z^{k},\label{eq:2.1}\end{equation}
with the parameters satisfying\begin{equation}
0\le b_{1},...,b_{s}<1\quad l>0.\label{eq:2.2}\end{equation}
It is clear that the function\begin{equation}
_{r}\phi_{s}\left(\begin{array}{c|c}
\begin{array}{c}
a_{1},\dotsc,a_{r}\\
b_{1},\dots,b_{s}\end{array} & q,-zq^{(s+1-r)/2}\end{array}\right)\label{eq:2.3}\end{equation}
 is of the form \eqref{eq:2.1} with \begin{equation}
l=\frac{s+1-r}{2}.\label{eq:2.4}\end{equation}
 Assume \begin{equation}
c:=c(a_{1},...,a_{r};b_{1},...,b_{s};q)=\frac{(-|a_{1}|,...,-|a_{r}|;q)_{\infty}}{(b_{1},...,b_{s};q)_{\infty}},\label{eq:2.5}\end{equation}
then,\begin{equation}
|f(z)|\le c\sum_{k=0}^{\infty}\frac{q^{lk^{2}}|z|^{k}}{(q;q)_{k}}.\label{eq:2.6}\end{equation}
We rewrite \eqref{eq:2.6} as\begin{equation}
c\sum_{k=0}^{\infty}\frac{q^{kl}(q^{l(k-1)}|z|)^{k}}{(q;q)_{k}}.\label{eq:2.7}\end{equation}
 For $k=0,1,...$, the terms\begin{equation}
\left(q^{l(k-1)}|z|\right)^{k}\label{eq:2.8}\end{equation}
 are bounded by\begin{equation}
\sqrt[4]{|z|^{2}q^{-l}}\,\exp\left(-\frac{\log^{2}|z|}{4l\log q}\right),\label{eq:2.9}\end{equation}
 thus,\begin{equation}
|f(z)|\le c\sqrt[4]{|z|^{2}q^{-l}}\,\exp\left(-\frac{\log^{2}|z|}{4l\log q}\right)\sum_{k=0}^{\infty}\frac{q^{kl}}{(q;q)_{k}},\label{eq:2.10}\end{equation}
 or\begin{equation}
|f(z)|\le\frac{c}{(q^{l};q)_{\infty}}\sqrt[4]{|z|^{2}q^{-l}}\,\exp\left(-\frac{\log^{2}|z|}{4l\log q}\right).\label{eq:2.11}\end{equation}
\[
\]

\begin{thm}
Assume that the entire function $f(z)$ is defined as in \eqref{eq:2.1}
and \eqref{eq:2.2}, then,\begin{equation}
|f(z)|\le\frac{(-|a_{1}|,...,-|a_{r}|;q)_{\infty}}{(b_{1},...,b_{s},q^{l};q)_{\infty}}\left\{ \sqrt[4]{|z|^{2}q^{-l}}\,\exp\left(-\frac{\log^{2}|z|}{4l\log q}\right)\right\} \label{eq:2.19}\end{equation}
 for any nonzero complex number $z$. In particular, for \begin{equation}
{0\le b}_{1},...,b_{s}<1,\quad0<s+1-r,\label{eq:2.20}\end{equation}
 we have\begin{gather}
\left|_{r}\phi_{s}\left(\begin{array}{c|c}
\begin{array}{c}
a_{1},\dotsc,a_{r}\\
b_{1},\dots,b_{s}\end{array} & q,z\end{array}\right)\right|\le\frac{(-|a_{1}|,...,-|a_{r}|;q)_{\infty}}{(b_{1},...,b_{s},q^{(s+1-r)/2};q)_{\infty}}\label{eq:2.21}\\
\times\left\{ \sqrt[4]{|z|^{2}q^{3(r-s-1)/2}}\,\exp\left(\frac{\log^{2}\left[|z|q^{(r-s-1)/2}\right]}{2(r-s-1)\log q}\right)\right\} \nonumber \end{gather}
 for any nonzero complex number.
\end{thm}
The same proof could be employed to show the following results.

\begin{thm}
Let $f(z)$ be a meromorphic function with the following Laurent expansion\begin{equation}
f(z)=\sum_{k=-\infty}^{\infty}a_{k}(z-a)^{k}\label{eq:2.22}\end{equation}
 on the whole complex plane except at $z=a$. If there exist two positive
numbers $\alpha>0$ and $0<q<1$ such that\begin{equation}
c=\sum_{k=-\infty}^{\infty}|a_{k}|q^{-|k|^{\alpha+1}}<\infty,\label{eq:2.23}\end{equation}
 then\begin{equation}
|f(z)|\le c\exp\left\{ \beta\left[\left|\log|z-a|\right|\right]^{\gamma}\right\} ,\label{eq:2.24}\end{equation}
 for $z\neq a$, where 

\begin{equation}
\beta=\frac{\alpha}{(\alpha+1)^{1+1/\alpha}\log^{1/\alpha}q^{-1}},\label{eq:2.25}\end{equation}
 and\begin{equation}
\gamma=\frac{\alpha+1}{\alpha}.\label{eq:2.26}\end{equation}
 In the particular, the theta function $\Theta(z|q)$ has an upper
bound\begin{equation}
|\Theta(z|q)|\le c\exp\left(\frac{\left|\log|z|\right|^{2}}{\log q^{-1}}\right)^{1/\alpha}\label{eq:2.28}\end{equation}
 for any $0<\alpha<1$, where\begin{equation}
c=\sum_{k=-\infty}^{\infty}q^{k^{2}-|k|^{1+\alpha}}.\label{eq:2.29}\end{equation}
 
\end{thm}
\begin{proof}
Clearly, we have\begin{equation}
|f(z)|\le\sum_{k=-\infty}^{\infty}|a_{k}|q^{-|k|^{\alpha+1}}\left(q^{|k|^{\alpha+1}}|z-a|^{k}\right).\label{eq:2.27}\end{equation}
For $k\in\mathbb{Z}$ and $z\neq a$, then, all the terms\begin{equation}
q^{|k|^{\alpha+1}}|z-a|^{k}\label{eq:2.28}\end{equation}
 are bounded by\begin{equation}
\exp\frac{\alpha\left[\left|\log|z-a|\right|\right]^{1+1/\alpha}}{(\alpha+1)^{1+1/\alpha}\log^{1/\alpha}q^{-1}},\label{eq:2.29}\end{equation}
 and the inequality \eqref{eq:2.24} follows. 
\end{proof}
\[
\]

\end{document}